\documentclass{amsart}
\usepackage{latexsym,amsmath, amssymb, amsfonts, amscd, amsthm, mathrsfs}

\renewcommand{\bar}{\overline}
\newcommand{\pC}{{{}^pC}}
\newcommand{\tC}{{{}^tC}}
\newcommand{\beC}{{{}^\beta C}}
\newcommand{\pM}{{{}^p\bM}}
\newcommand{\tM}{{{}^t\bM}}

\newcommand{\wH}{{{}^w H}}

\newcommand{\Ae}{{\mathbb A}}

\newcommand{\Ce}{{\mathbb C}}
\newcommand{\Fe}{{\mathbb F}}

\newcommand{\Qe}{{\mathbb Q}}
\renewcommand{\Re}{{\mathbb R}}
\newcommand{\Ze}{{\mathbb Z}}
\newcommand{\Cp}{{\mathbb C}_p}
\newcommand{\Fp}{{\mathbb F}_p}
\newcommand{\Qp}{{\mathbb Q}_p}
\newcommand{\Zp}{{\mathbb Z}_p}

\newcommand{\bA}{{\bf A}}

\newcommand{\bG}{{\bf G}}
\newcommand{\bK}{{\bf K}}
\newcommand{\bM}{{\bf M}}

\newcommand{\bP}{{\bf P}}

\newcommand{\br}{{\bf r}}

\newcommand{\cH}{ {\mathcal H} }
\newcommand{\cL}{ {\mathcal L} }

\newcommand{\cN}{ {\mathcal N} }
\newcommand{\cO}{ {\mathcal O} }
\newcommand{\cP}{ {\mathcal P} }

\newcommand{\cU}{ {\mathcal U} }
\newcommand{\cV}{ {\mathcal V} }

\newcommand{\M}{{\bf M}}
\newcommand{\Gal}{{\rm Gal}}

\newcommand{\barP}{\bar{P}}
\newcommand{\hP}{{\hat{P}}}
\newcommand{\hf}{{\hat{f}}}
\newcommand{\Tr}{  {\rm Tr}}
\newcommand{\ord}{ {\rm ord}}

\newcommand{\NP }{ {\rm  NP}}     %Newton polygon
\newcommand{\HP }{ {\rm  HP}}     %Hodge polygon
\newcommand{\sgn}{ {\rm sgn}}

\newcommand{\lcm}{{\rm lcm}}

\theoremstyle{plain}
\newtheorem{theorem}{Theorem}[section]
\newtheorem{proposition}[theorem]{Proposition}
\newtheorem{lemma}[theorem]{Lemma}

\newtheorem{corollary}[theorem]{Corollary}

\theoremstyle{definition}

\newtheorem{remark}[theorem]{Remark}
\theoremstyle{remark}
\newtheorem{acknowledgments}{Acknowledgments}

\begin{document}

\title
[Zeta functions of totally ramified $p$-covers]
{Zeta functions of totally ramified $p$-covers of the projective line}

 \subjclass[2000]{11,14}

\author{Hanfeng Li}
\address{
Department of mathematics, University of Toronto, Toronto, ON M5S 3G3,
CANADA
}
\email{hli@fields.toronto.edu}
\author{Hui June Zhu}
\address{
Department of maths and stats,
McMaster University,
Hamilton, ON L8S 4K1, CANADA
}
\email{zhu@cal.berkeley.edu}
\date{\today}

\begin{abstract}
In this paper we prove that there exists a Zariski dense open subset
$\cU$ defined over the rationals $\Qe$ in the space of
all one-variable rational functions with arbitrary $\ell$ poles of 
prescribed orders, such that for every geometric point $f$
in $\cU(\bar\Qe)$, the $L$-function of the exponential sum of $f$ at a
prime $p$ has Newton polygon approaching the Hodge polygon as $p$
approaches infinity. As an application to algebraic geometry, we prove
that the $p$-adic Newton polygon of the zeta function of a $p$-cover
of the projective line totally ramified at arbitrary $\ell$ points
of prescribed orders has an asymptotic generic lower bound.

\end{abstract}

\maketitle

\section{Introduction}

This paper investigates the asymptotics of the zeta functions of
$p$-covers of the projective line which are totally (wildly) ramified
at arbitrary $\ell$ points. Our approach is via Dwork's method on
one-variable exponential sums.

Throughout this paper we fix positive integers $\ell,
d_1,\ldots,d_\ell$, and let $d:=\sum_{j=1}^{\ell}d_j+\ell-2$. For
simplicity we assume $d\geq 2$ if $\ell=1$. Let $P_1=\infty$,
$P_2=0$, $P_3,\ldots,P_\ell$ be fixed poles in the projective line
over $\bar\Qe$ of orders
$d_1,\ldots,d_\ell$, respectively. Let $f$ be a one-variable
function over $\bar\Qe$ with these prescribed $\ell$ poles. It can
be written in a unique form of partial fractions
\cite[Introduction]{Zhu:1}):
\begin{eqnarray}\label{E:f}
f &=& \sum_{i=1}^{d_1}a_{1,i}x^i +
\sum_{j=2}^{\ell}\sum_{i=1}^{d_j}a_{ji}(x-P_j)^{-i}
\end{eqnarray}
with
$a_{ji}\in\bar\Qe$.
(Remark: we have assumed that $f$ has a vanishing constant term because
this does not affect the $p$-adic Newton polygons of $f$.)
Let $\Ae$ be the space of $a_{ji}$'s with
$\prod_{j=1}^{\ell}a_{j,d_j}\neq 0$.
It is an affine
($\sum_{j=1}^{\ell}d_j$)-space over $\Qe$. Let the Hodge polygon of $\Ae$,
denoted by $\HP(\Ae)$, be the lower convex graph of
the piecewise-linear function defined on the interval $[0,d]$ passing
through the two endpoints $(0,0)$ and $(d,d/2)$ and assuming every
slope in the list below of (horizontal) length $1$:
\begin{equation*}
\overbrace{0,\ldots,0}^{\ell-1}; \overbrace{1,\ldots,1}^{\ell-1};
\overbrace{\frac{1}{d_1}, \cdots, \frac{d_1-1}{d_1}}^{d_1-1};
\overbrace{\frac{1}{d_2}, \cdots, \frac{d_2-1}{d_2}}^{d_2-1};
\ldots\ldots; \overbrace{\frac{1}{d_\ell}, \cdots,
\frac{d_\ell-1}{d_\ell}}^{d_\ell-1}.
\end{equation*}
A non-smooth point on a polygon (as the graph of a piece-wise
linear function) is called a vertex. We remark that the classical
and geometrical `Hodge polygon' for any curve (including
Artin-Schreier curve as a special case) is the one with end points
$(0,0)$ and $(d,d/2)$ and one vertex at $(d/2,0)$. So the Hodge
polygon in our paper is different from the classical Hodge
polygon. We anticipate a $p$-adic arithmetic interpretation of our
Hodge polygon, but it remains an open question.

In \cite{Zhu:1} it is shown that in the case $\ell=1$ there is a
Zariski dense open subset $\cU$ defined over $\Qe$ such that every
geometric closed point $f$ in $\cU(\bar\Qe)$ has $p$-adic Newton polygon
approaching the Hodge polygon as $p$ approaches $\infty$. Wan
has proposed conjectures regarding multivariable exponential sums,
including the above as a special case (see \cite[Conjecture
1.15]{Wan:2}). This series of study traces back at least to
Katz \cite[Introduction]{Katz:1}, where Katz proposed to
study exponential sums in families instead of examining one
at a time. He systematically studied families of multivariable
Kloosterman sum in \cite{Katz:1}.

Let $\Qe_f$ be the extension field of $\Qe$ generated by coefficients 
$a_{ji}$'s and
poles $P_1,\ldots, P_\ell$ of $f$. For every prime number $p$ we
fix an embedding $\bar\Qe\hookrightarrow \bar\Qp$ once and for
all. This fixes a place $\cP$ in $\Qe_f$ lying over $p$ of residue
degree $a$ for some positive integer $a$. As usual, we let
$E(x)=\exp(\sum_{i=0}^{\infty} x^{p^i}/p^i)$ be the $p$-adic
Artin-Hasse exponential function. Let $\gamma$ be a root of the
$p$-adic $\log E(x)$ with $\ord_p(\gamma)=\frac{1}{p-1}$. Then
$E(\gamma)$ is a primitive $p$-th root of unity and we set
$\zeta_p:=E(\gamma)$. Let $\Fp$ be the
prime field of $p$ elements. Let $\Fe_q$ be a finite field of
$p^a$ elements. For $k\geq 1$, let $\psi_k: \Fe_{q^k} \rightarrow
\Qe(\zeta_p)^\times$ be a nontrivial additive character of
$\Fe_{q^k}$. Henceforth we fix $\psi_k(\cdot) =
\zeta_p^{\Tr_{\Fe_{q^k}/\Fe_p}(\cdot)}$. Let $\prod_{j=1}^{\ell}d_j$,
and all poles and leading coefficients $a_{j,d_j} $ of $f$ be
$p$-adic units. Let all coefficients $a_{j,i}$ of $f$ are $p$-adically
integral. (These are satisfied when $p$ is large enough.) Let
$S_k(f\bmod\cP) = \sum_x
\psi_k(f(x)\bmod\cP) $ where the sum ranges over all $x$ in
$\Fe_{q^k}\backslash \{\barP_1,\ldots,\barP_\ell\}$ (where
$\barP_j$ are
reductions of $P_j\bmod \cP$).
The $L$-function of $f$ at $p$ is defined as
\begin{eqnarray*}
L(f\bmod \cP;T) &=& \exp(\sum_{k=1}^{\infty}S_k(f\bmod\cP)T^k/k).
\end{eqnarray*}
This function lies in $\Ze[\zeta_p][T]$ of degree $d$. It is
independent of the choice of $\cP$ (that is, the embedding of
$\bar\Qe \hookrightarrow \bar\Qp$) for $p$ large enough, but we
remark that its Newton polygon is independent of the choice
of $\cP$ for all $p$ (see \cite[Section 1]{Zhu:26}). One notes
immediately that for every prime $p$ (coprime to the leading
coefficients, the poles and their orders) we have a map
$\NP_p(\cdot)$ which sends every $p$-adic integral point $f$
of $\Ae(\bar\Zp)$ to the Newton polygon $\NP_p(f)$ of the
$L$-function of exponential sums of $f$ at $p$. Given any
$f\in\Ae(\bar\Qe)$, we have for $p$ large enough that
$f\in\Ae(\bar\Zp)$ and hence we obtain the Newton polygon
$\NP_p(f)$ of $f$ at $p$. Presently it is known that $\NP_p(f)$
lies over $\HP(\Ae)$ for every $p$. These two polygons do not
always coincide. (See \cite[Introduction]{Zhu:26}.) Some
investigation on first slopes suggests the behavior is exceptional
if $p$ is small (see \cite[Introduction]{SZ:1}). There has been
intensive investigation on how the (Archimedean) distance between
$\NP_p(f)$ and $\HP(\Ae)$ on the real plane $\Re^2$ varies when
$p$ approaches infinity. Inspired by Wan's conjecture
\cite[Conjecture 1.15]{Wan:2} (proved in \cite{Zhu:1} for the
one-variable polynomial case), we believe that ``almost all"
points $f$ in $\Ae(\bar\Qe)$ satisfies $\lim_{p\rightarrow
\infty}\NP_p(f)=\HP(\Ae)$. Our main result is the following.

\begin{theorem}\label{T:main}
Let $\Ae$ be the coefficients space $\{a_{ji}\}$ of the $f$'s as in
(\ref{E:f}). There is a Zariski dense open subset
$\cU$ defined over $\Qe$ in
$\Ae$ such that for every geometric closed point $f$ in
$\cU(\bar\Qe)$ one has $f\in\cU(\bar\Zp)$ for $p$ large enough
(only depending on $f$), and
$$\lim_{p\rightarrow \infty}\NP_p(f) = \HP(\Ae).$$
\end{theorem}

The two polygons $\NP_p(f)$ and $\HP(\Ae)$ coincide if and only if
$p\equiv 1\bmod \lcm(d_j)$ (see \cite[Theorem 1.1]{Zhu:26}).
The case $\ell=1$ is known from \cite[Theorem 1.1]{Zhu:1}.
For $p\not\equiv 1\bmod \lcm(d_j)$,
the point
$f=x^{d_1}+\sum_{j=1}^{\ell} (x-P_j)^{-d_j}$
does not lie in $\cU$. This means $\cU$ is always a
proper subset of $\Ae$.

For any $\overline{f}\in\Ae(\Fe_q)$
and the (generalized)
Artin-Schreier curve
$C_{\overline{f}}: y^p-y = \overline{f}$,
let $\NP(C_{\overline{f}};\Fe_q)$ be the usual $p$-adic Newton polygon
of the numerator of the zeta function of $C_{\overline{f}}/\Fe_q$.
If it is shrunk by a factor of $1/(p-1)$ vertically and
horizontally, we denote it by
$\frac{\NP(C_{\overline{f}};\Fe_q)}{p-1}$.

\begin{corollary}\label{C:1}
Let notations be as in Theorem \ref{T:main} and the above.
For any $\overline{f}\in \Ae(\Fe_q)$
we have $\frac{\NP(\overline{f};\Fe_q)}{p-1}$ lies over
$\HP(\Ae)$ with the same endpoints, and they coincide if and only if
$p\equiv 1\bmod (\lcm(d_j))$. Moreover, there is a
Zariski dense open subset $\cU$
defined over $\Qe$ in $\Ae$ such that for every geometric closed point
$f$ in $\cU(\overline{\Qe})$ one has $f\in \cU(\overline{\Zp})$
for $p$ large enough (only depending on $f$), and
$$\lim_{p\rightarrow\infty}\frac{\NP(C_{\overline{f}};\Fe_q)}{p-1}
= \HP(\Ae).$$
\end{corollary}

\begin{proof}
This follows from the theorem above and a similar argument as
the proof of Corollary 1.3 in \cite{Zhu:26}, which we shall omit here.
\end{proof}

\begin{remark}\label{R:1}
(1) The result in Theorem \ref{T:main} and Corollary \ref{C:1} does not depend
on where those $\ell$ poles are
(as long as they are distinct).

(2) By Deuring-Shafarevic formula
(see for instance \cite[Corollary 1.5]{Crew:1}),
one knows that
$\NP_p(f)$ always has slope-$0$ segment precisely
of horizontal length $\ell-1$.
By symmetry it also has slope-$1$ segment of
the same length. See Remark 1.4 of \cite{Zhu:26}.
\end{remark}

Plan of the paper is as follows: section 2 introduces sheaves of
(infinite dimensional) $\varphi$-modules over some affinoid algebra arising
from one-variable exponential sums. We consider two Frobenius maps
$\alpha_1$ and $\alpha_a$. Section 4 is the main technical part, where
major combinatorics of this paper is done.  After
working out several combinatorial observations
we are able to reduce our problem to an
analog of the one-variable polynomial case as that in
\cite{Zhu:1}. Now back to Section 3 we improve the key lemma 3.5 of
\cite{Zhu:1} to make the generic Fredholm polynomial straightforward
to compute. Section 5 uses $p$-adic Banach theory to give a new
transformation theorem from $\alpha_1$ to $\alpha_a$ for any $a\geq
1$. This approach is very different from \cite{Wan:1} or
\cite{Zhu:2}. It shreds some new light on $p$-adic
approximations of $L$-functions of exponential sums and we believe
that it will find more application in the future.
Finally at the end of section 5 we prove our main result
Theorem \ref{T:main}.

\begin{acknowledgments}
Zhu's research was partially supported by an NSERC Discovery grant
and the Harvard University. She thanks Laurent Berger and the
Harvard mathematics department for hospitality during her visit in
2003. The authors also thanks the referee for comments.
\end{acknowledgments}

\section{Sheaves of $\varphi$-modules
over affinoid algebra}

The purpose of this section is to generalize the trace formula (see
\cite[Section 2]{Zhu:2}) for an exponential sum to that for families
of exponential sums.
See \cite{BGR:1} for fundamentals in rigid geometry and see
\cite{Berthelot} for an excellent setup for rigid cohomology
related to $p$-adic Dwork theory.

Let $\cO_1:=\Zp[\zeta_p]$ and $\Omega_1:=\Qp(\zeta_p)$. Fix a positive
integer $a$. Let $\Omega_a$ be the unramified extension of $\Omega_1$
of degree $a$ and $\cO_a$ its ring of integers.  Let
$\hP_1,\ldots,\hP_\ell$ in $\cO_a^\times$ be Teichm\"uller lifts of
$\bar P_1,\ldots,\bar P_\ell$ in $\Fe_{p^a}$. Similarly let $A_{j,i}$
be that of $\bar{a}_{j,i}$ and let $\vec{A}$ denote the sequence of
$A_{j,i}$ (we remark that for most part of the paper $\vec{A}$ will be
treated as a variable).  Let $\tau$ be the lift
of Frobenius to $\Omega_a$ which fixes
$\Omega_1$. Then $\tau(A_{ji})=A_{ji}^{p}$.  Let $1\leq j\leq
\ell$. Pick a root $\gamma^{1/d_j}$ of $\gamma$ in $\bar\Qp$ (or in
$\bar\Zp$, all the same) for the rest of the paper, and denote
$\Omega_1':=\Omega_1(\gamma^{1/d_1},\ldots,\gamma^{1/d_j})$. Let
$\cO_1'$ be its ring of integers. Let $\Omega_a':=\Omega_a\Omega_1'$
and let $\cO_a'$ be its ring of integers.  Then the affinoid algebra
$\cO_a'\langle\vec{A}\rangle$ (with $\vec{A}$ as variables) forms a
Banach algebra over $\cO_a'$ under the supremum norm.

Let $0<r<1$ and $r\in |\Omega_a'|_p$. Let $\bA_r$ be the affinoid
with $\ell$ deleted discs centering at $\hP_1,\ldots,\hP_\ell$
each of radius $r$ on the rigid projective line $\bP^1$ over $\Omega_a'$
(as defined in
\cite{Zhu:26}). The topology on $\bA_r$ is given by
the fundamental system of strict neighborhood $\bA_{r'}$ with
$r\leq r'<1$ and $r'\in |\Omega_a'|_p$.
Let $\bA$ be $\bA_r$ for some unspecified $r$
sufficiently close to $1^{-}$
(the precise bound on the size of $r$ was discussed in
\cite[Section 2]{Zhu:26}).
Let $\cH(\Omega_a')$ be the ring of rigid analytic functions on $\bA$
over $\Omega_a'$. Then it is a $p$-adic Banach space over $\Omega_a'$.
It consists of functions in one variable $X$ of the form
$\xi=\sum_{i=0}^{\infty}c_{1,i}X^i
+ \sum_{j=2}^{\ell}\sum_{i=1}^{\infty}c_{j,i}(X-\hP_j)^{-i}$
where $c_{j,i}\in\Omega_a'$
and $\forall j\geq 1, \lim_{i\rightarrow \infty} \frac{|c_{j,i}|_p}{r^i}=0$.
Its norm is defined as $||\xi||=\max_{j}(\sup_{i}\frac{|c_{j,i}|_p}{r^i})$.
(See \cite[Section 2.1]{Zhu:26}.)
  %Note that $\cH_\bA$
  %is functorial (see discussion in \cite{Colmez:1}).
Let $\cH(\Omega_a'\langle\vec{A}\rangle) := \cH(\Omega_a')
\hat\otimes_{\Omega_a'}\Omega_a'\langle\vec{A}\rangle$ where
$\hat\otimes$ means $p$-adic completion after the tensor product. It is a
$p$-adic Banach modules over $\Omega_a'\langle\vec{A}\rangle$ with
the natural norm on the tensor product of two Banach modules defined by
the followings.
For any $\sum v\otimes w \in \cH(\Omega_a')\otimes
\Omega_a\langle\vec{A}\rangle$ let
$||\sum v\otimes w||=
\inf(
\max_{i}(||v_i||\cdot||w_i ||)
)$,
where the $\inf$ ranges over all representatives $\sum_i v_i\otimes w_i$ with
$\sum v\otimes w = \sum_i v_i \otimes w_i$.
From the
$p$-adic Mittag-Leffler decomposition theorem derived in
\cite[Section 2.1]{Zhu:26}, we can generalize it to the decomposition of
$\Omega_a'\langle\vec{A}\rangle$ as a Banach
$\Omega_a'\langle\vec{A}\rangle$-module.
Write $X_1=X$ or $X_j = (X-\hP_j)^{-1}$ for $2\leq j\leq \ell$.
Let $Z_j=\gamma^{1/d_j}X_j$. Note that $\vec{b}_{\rm w}=\{1, Z^i_1, \ldots,
Z^i_{\ell}\}_{i\ge 1}$ is a formal basis of the Banach
$\Omega_a'\langle\vec{A}\rangle$-module
$\cH(\Omega_a'\langle\vec{A}\rangle)$, that is, every $v$ in
$\cH(\Omega_a'\langle\vec{A}\rangle)$ can be written uniquely as
an infinite sum of $c'_{j,i}Z_j^i$'s with $c'_{j,i}\in
\Omega_a'\langle\vec{A}\rangle$ and $\frac{|c'_{j,i}|_p}{r'^i}
\rightarrow 0$ as
$i\rightarrow \infty$, where $r'=r\;p^{-\frac{1}{(p-1)d_j}}$.
The Banach module
$\cH(\Omega_a'\langle\vec{A}\rangle)$ is orthonormalizable (even
though $\vec{b}_{\rm w}$ is not its orthonormal basis).

In this paper we extend the
$\tau$-action so that it acts on $\gamma^{\frac{1}{d_J}}$ trivially for any $J$.
Below we begin to construct the Frobenius operator $\alpha_1$ on
$\cH(\Omega_a'\langle\vec{A}\rangle)$.
Recall the $p$-adic Artin-Hasse exponential function
$E(X)$. Take expansion of $E(\gamma X)$ at $X$ one gets
$E(\gamma X) =\sum_{m=0}^{\infty}\lambda_m X^m$ for some
$\lambda_m\in\cO_1$. Clearly $\ord_p\lambda_m\geq
\frac{m}{p-1}$ for all $m\geq 0$. In particular, for $0\leq m\leq
p-1$ the equality holds and $\lambda_m = \frac{\gamma^m}{m!}$.
Let $ F_j(X_j):= \prod_{i=1}^{d_j} E(\gamma A_{j,i}X_j^i)$. Then
$$
F_j(X_j)
=\sum_{n=0}^{\infty}
F_{j,n}(A_{j,1},\ldots,A_{j,d_j})X_j^n,
$$
where $F_{j,n}:=0$ for $n<0$ and for $n\geq 0$
\begin{eqnarray}\label{E:F}
F_{j,n} &:=&
  \sum \lambda_{m_{1}}\cdots\lambda_{m_{d_j}}
  A_{j,1}^{m_1}\cdots A_{j,d_j}^{m_{d_j}},
\end{eqnarray}
where the sum ranges over all $m_1,\ldots,m_{d_j}\geq 0$ and
$\sum_{k=1}^{d_j}km_k = n$. It is clear that $F_{j,n}$ lies in $
\cO_1[A_{j,1},\ldots,A_{j,d_j}] $. One observes that $F_j(X_j)\in
\cO_1\langle A_{j,1},\ldots,A_{j,d_j}\rangle \langle X_j\rangle$,
the affinoid algebra in one variable $X_j$ (actually it lies in
$\cO_1[A_{j,1},\ldots,A_{j,d_j}]\langle X_j\rangle)$. Taking
product over $j=1,\ldots,\ell$, we have that $ F(X):=
\prod_{j=1}^{\ell} F_j(X_j) $ lies in
$\cH(\cO_a\langle\vec{A}\rangle)$. Let $\tau_*^{a-1}$ be the
push-forward map of $\tau^{a-1}$, that is, for any function $\xi$,
$\tau_*^{a-1}(\xi) = \tau^{a-1}\circ \xi \circ \tau$. For example,
$\tau_*^{a-1}(B/(X-\hP^p)) = \tau^{a-1}(B)/(X-\hP)$ for any
$B\in\cO_1\langle \vec{A}\rangle$ and $\hP$ a Teichm\"uller lift
of some $\barP$. Let $U_p$ be the Dwork operator and let $F(X)$
denote the multiplication map by $F(X)$, as defined in
\cite[Section 2]{Zhu:26}. Let $\alpha_1:=\tau_*^{a-1} \circ
U_p\circ F(X)$ denote the composition map.  Then $\alpha_1$ is a
$\tau^{a-1}$-linear endomorphism of $\cH(\Omega_a'\langle\vec{A}
\rangle)$ as a Banach $\Omega_a'\langle\vec{A}\rangle$-module.

Let $S$ be the affinoid over $\Omega_a'$ with affinoid algebra
$\Omega_a'\langle\vec{A}\rangle$.  If $\cL$ is a sheaf of $p$-adic
Banach $\Omega_a'\langle\vec{A}\rangle$-module (with formal basis) and
$\alpha_1$ is the Frobenius map which is $\tau^{a-1}$-linear with
respect to $\Omega_a'\langle\vec{A}\rangle$, then we call the pair
$(\cL,\alpha_1)$ a {\em sheaf of $\varphi$-module of infinite
rank}. Note that the pair
$(\cH(\Omega_a'\langle\vec{A}\rangle),\alpha_1)$ can be considered as
sections of a sheaf of $\Omega_a'\langle\vec{A}\rangle$-module of
infinite rank on $\bA$. This is intimately related to Wan's {\em
nuclear $\sigma$-module of infinite rank} (see \cite{Wan:3}) if
replacing his $\sigma$ by our $\tau^{a-1}$. Wan has defined
$L$-functions of nuclear $\sigma$-modules and he also showed that it
is $p$-adic meromorphic on the closed unit disc (see Wan's papers
\cite{Wan:3,Wan:4} which proved Dwork's conjecture).
Finally we define $\alpha_a:=\alpha_1^a$.

Recall that $\alpha_1$ is a $\tau^{a-1}$-linear (with respect to
$\Omega_a'\langle\vec{A}\rangle$) completely continuous endomorphism on
the $p$-adic Banach module $\cH(\Omega_a'\langle\vec{A}\rangle)$ over
$\Omega_a'\langle \vec{A}\rangle$. Let $1\leq J_1,J\leq \ell$. Write
$(\alpha_1 Z_J^i)_{\hP_{J_1}} = \sum_{n=0}^{\infty} (\tau^{a-1}
C_{J_1,J}^{n,i}) Z_{J_1}^n$ for some $C_{J_1,J}^{n,i}$ in
$\Omega_a'\langle\vec{A}\rangle$. The matrix of $\alpha_1$,
consisting of all these $\tau^{a-1}C_{J_1,J}^{\star,\star}$'s, is
a nuclear matrix (see section 5). This matrix is the subject of
the next section.
Below we extend Dwork, Monsky and Reich's trace formula
to families of one-variable exponential sums.

\begin{theorem}\label{T:main2}
Let $\overline{f}=\sum_{i=1}^{d_1}\overline{a}_{1,i}x^i
+\sum_{j=2}^{\ell}\sum_{i=1}^{d_i}\overline{a}_{j,i}(x-P_j)^{-i}
\in \Ae (\Fe_q)$ and let
$\hf$ be its Teichm\"uller lift
with coefficient $\overline{a}_{ji}$ being lifted to $A_{ji}$.
Let $\cH(\Omega_a'\langle\vec{A}\rangle)_r$
be the Banach module $\cH(\Omega_a'\langle\vec{A}\rangle)$
for some suitably chosen $0<r<1$ with $r\in|\Omega_a'|_p$
close enough to $1^{-}$.
Then
\begin{eqnarray*}
L(\overline{f}/\Fe_q; T)
&=& \frac{\det(1-T\alpha_a|\cH(\Omega_a'\langle \vec{A}\rangle))}
{\det(1-Tq\alpha_a|\cH(\Omega_a'\langle\vec{A}\rangle))}
\end{eqnarray*}
lies in $\cO_a\langle \vec{A}\rangle [T]$ as a polynomial of
degree $d$ in $T$. Its Teichm\"uller specialization of $\vec{A}$ in
$\cO_a$ lies in $\Ze[\zeta_p][T]$.
\end{theorem}
\begin{proof}
The proof is similar to that of \cite[Lemma 2.7]{Zhu:2}. Let
$\cH^\dagger:=\bigcup_{0<r<1}\cH(\Omega_a'\langle\vec{A}\rangle)_r$.
Then it is the Monsky-Washnitzer dagger space. Then $\alpha_a$ is
a completely continuous endomorphism on $\cH^\dagger$ and the
determinant $\det(1-T \alpha_a|\cH^\dagger) =
\det(1-T\alpha_a|\cH(\Omega_a'\langle\vec{A}\rangle)_r)$ for any
$r$ within suitable range in $(0,1)$ is independent of $r$.
Finally one knows that the coefficients are all integral so lies
in $\cO_a$ and coefficient of $T^m$ vanishes for all $m>d$.
We omit details of the proof.
\end{proof}

\section{Explicit approximation of the Frobenius matrix}

This section uses some standard techniques in $p$-adic
approximation and it is very technical. The readers are
recommended to skip it at first reading and continue at
the next section.

\subsection{The nuclear matrix}\label{S:3.1}

Let notation be as in the previous section. Assign $\phi(1)=0$.
Let $\phi(Z^n_j)=\frac{n}{d_j}$ for $j\le 2$ or $\frac{n-1}{d_j}$
for $j\ge 3$. Order the elements in $\vec{b}_{\rm w}$ as
$e_1, e_2, \cdots$ such that $\phi(e_1)\le \phi(e_2)\le \cdots$.
Consider the infinite matrix representing the endomorphism
$\alpha_1$ of the $\Omega_a'\langle\vec{A}\rangle$-module
$\cH(\Omega_a'\langle\vec{A}\rangle)$ with respect to the basis
$\vec{b}_{\rm w}$. This matrix can be written as $\tau^{a-1}\bM$,
where each entry is $\tau^{a-1}C_{J_1,J}^{n,i}$ for  $1\leq
J_1,J\leq \ell$.

Our goal of this section is to collect delicate information about
entries of the matrix $\bM$. Recall the is polynomial $F_{J,n_J}$ in
$\cO_1[\vec{A}]$ as in (\ref{E:F}), which we have already built up
some satisfying knowledge. Below we will express
$C_{J_1,J}^{\star,\star}$ as a polynomial expression in these
$F_{J_1,n_{J_1}}$'s. In this paper the {\em formal expansion} of
$C_{J_1,J}^{\star,\star}$ will always mean the formal sum in
$\cO_a'[\vec{A}]$ by the composition of (\ref{E:CJJ}) and the
formula in Lemma \ref{L:Hjj}.

For $n,i\geq 1$, and if $J=1$ or $J_1=1$
then for $i\geq 0$ or for $n\geq 0$ respectively one
has
\begin{equation}\label{E:CJJ}
C_{J_1,J}^{n,i} =
\left\{
\begin{array}{ll}
\gamma^{\frac{i}{d_J}-\frac{n}{d_{J_1}}}
H_{J_1,J}^{np,i}
 & J_1=1,2\\
\gamma^{\frac{i}{d_J}-\frac{n}{d_{J_1}}} \sum_{m=n}^{np}
C^{n,m}\hP_{J_1}^{np-m}
                H_{J_1,J}^{m,i}

 & J_1\geq 3
\end{array}
\right.
\end{equation}
where $C^{\star,\star}\in \Zp$ is defined in \cite[Lemma
3.1]{Zhu:26} and $H_{J_1,J}^{\star,\star} \in
\cO_a\langle\vec{A}\rangle$ is formulated in Lemma \ref{L:Hjj}
below. Indeed, we recall that $C^{n,m}$ is actually a rational
integer and it only depends on $n$,$m$ and $p$.

\begin{lemma}\label{L:Hjj}
Let $\vec{n}:=(n_1,\ldots,n_\ell)\in\Ze_{\geq 0}^\ell$.

(1) For $i,n\geq 0$, then $H_{1,J}^{n,i}$ is equal to
\begin{eqnarray*}
%\label{E:P1}
\sum
\left( F_{1,n_1}\cdot
\left(\sum_{
\stackrel{J\neq 1}{0\leq m_J\leq n_J}}
F_{J,m_J}\binom{n_J+i-1}{m_J+i-1}\hP_J^{n_J-m_J}
\right)\right.\\
\cdot \left.\prod_{j\neq 1,J}
  \left(\sum_{m_j=0}^{n_j} F_{j,m_j}\binom{n_j-1}{m_j-1}\hP_j^{n_j-m_j}
  \right)
  \right),
\end{eqnarray*}
where the sum ranges over all $\vec{n}\in \Ze_{\geq 0}^\ell$ such that
$n=n_1\pm i -\sum_{j=2}^{\ell}n_j$ and the $+$ or $-$ depends on
$J=1$ or $J\neq 1$, respectively.

(2) For $J_1,J\neq 1$, one has that $H_{J_1,J}^{n,i}$ is equal to
\begin{eqnarray*}
\sum
\left(
F_{J_1,n_{J_1}}\cdot
\left(\sum_{\stackrel{J\neq J_1}{m_J\geq 0}}
F_{J,m_J}(-1)^{m_J+i}
  \binom{n_J+m_J+i-1}{m_J+i-1}(\hP_J-\hP_{J_1})^{-(n_J+m_J+i)}
\right)
\right.
\\
\cdot
\left(\sum_{m_1=n_1}^{\infty}F_{1,m_1}\binom{m_1}{n_1}
\hP_{J_1}^{m_1-n_1}\right)\\
\left.
\prod_{j\neq 1,J_1,J}
\left(\sum_{m_j=0}^{\infty} F_{j,m_j}(-1)^{m_j}
  \binom{n_j+m_j-1}{m_j-1}(\hP_j-\hP_{J_1})^{-(n_j+m_j)}
  \right)
  \right)
\end{eqnarray*}
where the sum ranges over all $\vec{n}\in\Ze_{\geq 0}^\ell$ such that
$n=n_{J_1}+i-\sum_{j\neq J_1}n_j$ if $J=J_1$ and
$n=n_{J_1}-\sum_{j\neq J_1}n_j$ if $J\neq J_1$.

(3) For $J_1\neq 1$ and $J=1$ we have that $H_{J_1,J}^{n,i}$ is equal to
\begin{eqnarray*}
\sum
\left(
F_{J_1,n_{J_1}}\cdot
\left(\sum_{m_1=n_1-i}F_{1,m_1}\binom{m_1+i}{n_1}\hP_{J_1}^{m_1+i-n_1}
\right)
\right.
\\
\cdot
\left.
\prod_{j\neq J_1,1}
\left(\sum_{m_j=0}^{\infty} F_{j,m_j}(-1)^{m_j}
  \binom{n_j+m_j-1}{m_j-1}(\hP_j-\hP_{J_1})^{-(n_j+m_j)}
  \right)
  \right)
\end{eqnarray*}
where the sum ranges over all $\vec{n}\in\Ze_{\geq 0}^\ell$ such that
$n=n_{J_1}-\sum_{j\neq J_1}n_j$.
\end{lemma}

\begin{proof}
We shall use ``$\stackrel{\hP_j}{=}$" to mean expansion at $\hP_j$.
Clearly
for any $J_1$ one has
$F_{J_1}(X_{J_1})X_{J_1}^i
\stackrel{\hP_{J_1}}{=}
\sum_{n=0}^{\infty}F_{J_1,n}X^{n+i}$.

For $J\geq 2$ one has the expansion at
$\hP_1=\infty$:
\begin{eqnarray*}
F_J(X_J)X_J^i
 &=& \sum_{m=0}^{\infty}F_{J,m}(X^{-1}(1-\hP_JX^{-1})^{-1})^{m+i}\\
 &\stackrel{\hP_1}{=}& \sum_{m=0}^{\infty}F_{J,m}
  \sum_{k=m+i}^{\infty}\binom{k-1}{m+i-1}\hP_J^{k-(m+i)}X^{-k}\\
 &=& \sum_{n=0}^{\infty}(\sum_{m=0}^{n}F_{J,m}
     \binom{n+i-1}{m+i-1}\hP_J^{n-m})
     X^{-n-i}.
\end{eqnarray*}

For $J_1\neq 1$ and $J\neq 1,J_1$,
its expansion at $\hP_{J_1}$ is:
\begin{eqnarray*}
F_J(X_J)X_J^i
 &=& \sum_{m=0}^{\infty}F_{J,m}(X_{J_1}^{-1}-(\hP_J-\hP_{J_1}))^{-(m+i)}\\
 &\stackrel{\hP_{J_1}}{=}& \sum_{m=0}^{\infty}F_{J,m}
      (-1)^{m+i}\sum_{n=0}^{\infty}\binom{n+m+i-1}{m+i-1}
      (\hP_J-\hP_{J_1})^{-(n+m+i)}X_{J_1}^{-n}\\
 &=&\sum_{n=0}^{\infty}
       (\sum_{m=0}^{\infty}F_{J,m}(-1)^{m+i}\binom{n+m+i-1}{m+i-1}
        (\hP_J-\hP_{J_1})^{-(n+m+i)})X_{J_1}^{-n}.
\end{eqnarray*}
For $J_1\neq 1$ and $J=1$ then
one has
\begin{eqnarray*}
F_J(X)X_J^i
 &\stackrel{\hP_{J_1}}{=}
   &\sum_{n=0}^{\infty}
    \left(\sum_{m=n-i}^{\infty}F_{J,m}\binom{m+i}{n}\hP_{J_1}^{m+i-n}
    \right)X_{J_1}^{-n}.
\end{eqnarray*}

By $F(X)X_J^i=(F_J(X_J)X_J^i)\cdot\prod_{j\neq J}F_j(X_j)$,
and Key Computational Lemma of \cite{Zhu:26},
one can compute and obtain
$(F(X)X_J^i)_{\hP_{J_1}}$ for the case $J_1=1$ or $J_1\neq 1$
presented respectively in the two formulas in our assertion.
This proves the lemma.
\end{proof}

\begin{remark}
If we are dealing with the case of unique pole at $\infty$ then one
sees easily that $C_{1,1}^{\star,\star}$ lies in
$\cO_1'[\vec{A}]$. This greatly reduces the complexity of
situation.
\end{remark}

The following results were presented in \cite{Zhu:26}.
See Section 3 and in particular,
Theorem 3.7 of \cite{Zhu:26} for a proof.
We shall use $t_{J_1}$ to denote the
lower bound in Lemma \ref{L:F&H} c).

\begin{lemma}\label{L:F&H}
Let notation be as above.

(a) For all $J$ and $n_J$ we have
$\ord_p(F_{J,n_J})
  \geq \frac{\lceil\frac{n_J}{d_J}\rceil}{p-1}
  \geq \frac{n_J}{d_J(p-1)}.$

(b) For all $J_1,J$, and all $n,i$ we have
$\ord_p(H_{J_1,J}^{n,i})
  \geq \frac{n-i}{d_{J_1}(p-1)}.$

(c) For any $J_1$ and any $n$ we have
$
\ord_p (C_{J_1,\star}^{n,\star})
 \geq \frac{n}{d_{J_1}} \mbox{ or } \frac{n-1}{d_{J_1}}
$
depending on $J_1=1,2$ or $3\leq J_1\leq \ell$.
Moreover,
$
\ord_p(C_{J_1,J_1}^{n,i}) \geq
{\lceil\frac{np-i}{d_{J_1}}\rceil}/(p-1)
\mbox{ or } \lceil\frac{(n-1)p-(i-1)}{d_{J_1}}\rceil/(p-1)
$
depending on $J_1=1,2$ or $3\leq J_1\leq \ell$, respectively.
\end{lemma}

\subsection{Approximation by truncation}

{}From the previous subsection one has
noticed an unpleasant feature of
$C_{J_1,J}^{\star,\star}$ for the purpose of
approximation by $F_{J_1,n_{J_1}}$'s.
First, in the sum for (\ref{E:CJJ}) when $J_1\geq 3$,
the range of $m$ is too `large'.
Second, $H_{J_1,J}^{\star,\star}$
of Lemma \ref{L:Hjj} is generally an infinite sum
of $F_{J_1,n_{J_1}}$'s.
In this subsection
we will define an approximation
in terms of truncated finite sum
of $F_{J_1,n_{J_1}}$'s.
Below we prove two lemmas which will be used for approximation
in Lemma \ref{L:keylemma}.

For any integer $0<t\leq p$, let $\tC_{J_1,J}^{n,i}$
be the same as $C_{J_1,J}^{n,i}$
except for $J_1\geq 3$
its sum ranges over all $m$
in the sub-interval $[(n-1)p+1,(n-1)p+t]$.

\begin{lemma}\label{L:compare}
Let $3\leq J_1\leq \ell$,
$1\leq J\leq \ell$.
Let $n\leq d_{J_1}$ and $i\leq d_J$.

(1) For $p$ large enough, one has
\begin{eqnarray}
\ord_p(C_{J_1,J}^{n,i}-\pC_{J_1,J}^{n,i})
 &>&
\frac{n-1}{d_{J_1}}+\frac{d}{p-1}.
\end{eqnarray}

(2) There is a constant $\beta>0$ depending only on $d$ such that
for $t\geq \beta$ one has
\begin{eqnarray}
\ord_p(\pC_{J_1,J}^{n,i} - \tC_{J_1,J}^{n,i}) & > &
\frac{n-1}{d_{J_1}} +\frac{d}{p-1}.
\end{eqnarray}
\end{lemma}

\begin{proof}
(1) By \cite[Lemma 3.1]{Zhu:26}, one knows that for any $m\leq
(n-1)p$ one has $\ord_p(C^{n,m}) \geq 1$
and hence
$\ord_p(C_{J_1,J}^{n,i}) \geq 1+
(\frac{i}{d_J}-\frac{n}{d_{J_1}})\frac{1}{p-1}$.
For $n\leq d_{J_1}$ and for $p$ large
enough one has $1 + (\frac{i}{d_J}-\frac{n}{d_{J_1}})\frac{1}{p-1}
> \frac{n-1}{d_{J_1}} + \frac{d}{p-1}$.
Combining these two inequalities, one concludes.

(2)
We may assume $J_1\geq 3$. Then
for any $1\leq v\leq p$, by Lemma \ref{L:F&H},
\begin{eqnarray*}
\ord_p (H_{J_1,J}^{(n-1)p+v,i}) &\geq &
\frac{(n-1)p+v-i}{d_{J_1}(p-1)} >
(\frac{n}{d_{J_1}}-\frac{i}{d_J})\frac{1}{p-1}
    + \frac{n-1}{d_{J_1}} +\frac{d}{p-1},
\end{eqnarray*}
if $v\geq \beta$ for some $\beta>0$ only depending on $d$.
Therefore,
\begin{eqnarray*}
\ord_p(\pC_{J_1,J}^{n,i}-\tC_{J_1,J}^{n,i}) & > &
\frac{n-1}{d_{J_1}} + \frac{d}{p-1}.
\end{eqnarray*}
This finishes our proof.
\end{proof}

Fix $\beta$ for the rest of the paper. We will truncate the
infinite expansion of $H_{J_1,J}^{\star,\star}$. Let $w>0$ be any
integer. For $J_1=1,2$ let ${\wH}_{J_1,J}^{np,i}$ be the sub-sum
in $H_{J_1,J}^{np,i}$ where $\vec{n}=(n_1,\ldots,n_\ell)$ are such
that $n_{J_1}-np$ and $n_j$ lie the interval $[-w,w]$ for $j\neq
J_1$. Similarly, for $J_1\geq 3$ let ${\wH}_{J_1,J}^{m,i}$ be the
sub-sum of $H_{J_1,J}^{m,i}$ where $\vec{n}$ ranges over the
finite set of vectors $(n_1,\ldots,n_\ell)$ such that
$n_{J_1}-(n-1)p$ and $n_j$ lie in the interval $[-w,w]$ for $j\neq
J_1$. Consider $\beC_{J_1,J}^{n,i}$ as a polynomial expression in
$H_{J_1,J}^{\star,\star}$'s, then we set $ {{}^wK}_{J_1,J}^{n,i}
:= \beC_{J_1,J}^{n,i}({\wH}_{J_1,J}^{n,i}). $

\begin{lemma}\label{L:finite-term}
There is a constant $\alpha$
depending only on $d$ such that
\begin{eqnarray}
\ord_p( \beC_{J_1,J}^{n,i} - {{}^{\alpha}K}_{J_1,J}^{n,i})
 & > &
\frac{n-1}{d_{J_1}} +\frac{d}{p-1}.
\end{eqnarray}
\end{lemma}
\begin{proof}
This part is similar to Lemma \ref{L:compare} 2), so we omit its proof.
\end{proof}

\subsection{Minimal weight terms}

The {\em weight} of a monomial (with nonzero coefficient)
$(\prod_{j=1}^{\ell}\prod_{i=1}^{d_j} A_{j,i}^{k_{j,i}})$ in
$\cO_a'[\vec{A}]$ is defined as $\sum_{j=1}^{\ell}\sum_{i=1}^{d_j}
i k_{j,i}.$ For example, the weight of $A_{1,2}^aA_{1,3}^b$ is
equal to $2a+3b$. We will later utilize the simple observation
that every monomial in $F_{J,n_J}$ is of weight $n_J$.

We call those entries  with $J_1=J$ the diagonal one (or blocks).
As we have seen in Lemma \ref{L:Hjj}, the off-diagonal entries are
less manageable while the diagonal entries behave well in
principle. For any integer $0<t\leq p$, let
$\tM:=(\tC_{J_1,J}^{n,i})$ with respect to the basis arranged in
the same order as that for $\bM$. Consider the diagonal blocks,
consisting of $\pC_{J,J}^{\star,\star}$'s. Despite
$\pC_{J,J}^{\star,\star}$ lives in $\cO_a'\langle\vec{A}\rangle$,
its minimal weight terms live in
$\hP_J^\Ze\gamma^{\frac{i-n}{d_J}} \cO_1[\vec{A}_J]$.

\begin{lemma}\label{L:locate-min}
Let $p>d_j$ for all $j$. The  minimal weight monomials of
$\pC_{J,J}^{n,i}$ (with $J=1,2$) live in the term $\gamma^{\frac{i-n}{d_1}}
F_{1,np-i}$ where $d_1>n,i\geq 0$ unless $n=0$ and $i>0$.
For $J\geq 3$ and $n\geq 2$, the minimal weight monomials of
$\pC_{J,J}^{n,i}$ live in the term
$$
\gamma^{\frac{i-n}{d_J}}
C^{n,(n-1)p+1} \hP_J^{p-1}
F_{J,(n-1)p-(i-1)}
$$
where $d_J>n,i\geq 1$.
\end{lemma}

\begin{proof}
This follows from Lemma \ref{L:Hjj}. We omit its proof.
\end{proof}

Given a $k \times k$ matrix $M:=(m_{ij})_{1\leq i,j\leq k}$ with a
given formal expansion of $m_{ij} \in
\cO_a'\langle\vec{A}\rangle$, the {\em formal expansion} of $\det
M$ means the formal expansion as $\sum_{\sigma\in
S_k}\sgn(\sigma)\prod_{n=1}^{k} m_{ij}$ where the product is
expanded according to the given formal expansion $m_{ij}$. For
example, if $m_{ij}=C_{J_1,J}^{\star,\star}$ then its formal
expansion is given by composition of (\ref{E:CJJ}) and formulas in
Lemma \ref{L:Hjj}.

\begin{lemma}\label{L:minimal-weight}
Let notation be as above and let $p>d_j$ for all $j$. Then in the
formal expansion of $\det(\pM)^{[k]}$ in
$\cO_a'\langle\vec{A}\rangle$, all minimal weight terms are from
$\prod_{J=1}^{\ell}\det\pC_{J,J}^{n,i}$ (with $n,i\geq 1$ in a
suitable range for $J=1,2$ and with $n,i\geq 2$ for $J\geq 3$) of the
diagonal blocks.
\end{lemma}
\begin{proof}
We will show that picking an arbitrary entry on the diagonal block,
every off-diagonal entry on the same row has
strictly higher minimal weight among its monomials.

Let $\vec{A}_J$ stand for the vector $(A_{J,1},
\ldots,A_{J,d_J})$. As we have noticed earlier the polynomial
$F_{J,n_J}$ in $\cO_1[\vec{A}_J]$ has every monomial of equal
weight $n_J$ for any $J$. For simplicity we assume $n,i\geq 1$
here. Using data from Lemma \ref{L:Hjj}, we find all minimal
weight monomials in $H_{J_1,J}^{\star,\star}$'s illustrated below
by an arrow: $H_{1,1}^{np,i} \rightarrow F_{1,np-i}, H_{1,J\geq
2}^{np,i} \rightarrow  F_{1,np+i}, H_{2,2}^{np,i} \rightarrow
F_{2,np-i}, H_{2,J\neq 2}^{np,i} \rightarrow F_{2,np}$.
One also notes that for $J_1\geq 3$ one has
that $H_{J_1\geq 3,J}^{(n-1)p+1,i} \rightarrow  F_{J_1,(n-1)p-(i-1)}$
if  $J_1=J$, and $H_{J_1\geq 3,J}^{(n-1)p+1,i}
\rightarrow F_{J_1,(n-1)p+1}$ if $J_1\neq J$.
One notices from (\ref{E:CJJ}) and the above that
the minimal weight monomials of $\pC_{J_1,J}^{n,i}$
live in $H_{J_1,J}^{np,i}$ if $J_1=1,2$
and in $H_{J_1,J}^{(n-1)p+1,i}$ if $3\leq J_1\leq \ell.$

Recall that for $J=1$ the range for $i$ is $i\geq 0$. In all other
cases the range is $i\geq 1$. From the above we conclude our claim
in the beginning of the proof. Consequently, all minimal weight
monomials in the formal expansion of the determinant
$\det\bM^{[k]}$ come from the diagonal blocks. By
Lemma \ref{L:locate-min},
$C_{1,1}^{0,i}$ and $C_{J,J}^{1,i}$ (with $J\geq 3$)
both have their minimal weight equal to $0$ if $i=0$ and $>0$ if
$i>0$. Then it is not hard to conclude that the minimal weight
monomials of $\det (C_{1,1}^{n,i})_{n,i\geq 0}$
(resp. $\det (C_{J,J}^{n,i})_{n,i\geq 1}$)
are from $\det (C_{1,1}^{n,i})_{n,i\geq 1}$
(resp. $\det (C_{J,J}^{n,i})_{n,i\geq 2}$).
\end{proof}

For $1\leq J\leq \ell$, let $D_J^{[k]} := \det
(F_{J,ip-j})_{1\leq i,j \leq k} \in\cO_1[A_1,\ldots,A_d]$.

\begin{proposition}\label{P:weight}
Let $p>d_j$ for all $j$.
The minimal weight monomials of $\det((\pC_{J,J}^{i,j})_{1\leq
i,j\leq k}$ for $J=1,2$ (resp.
$\det((\pC_{J,J}^{i,j})_{2\leq
i,j\leq k}$ for $J\geq 3$
) lie in $D_J^{[k]}$ (resp. $D_J^{[k-1]}$).
Every monomial of $D_J^{[k]}$ (resp. $D_J^{[k-1]}$)
corresponds to a monomial in the formal expansion of $\det
((\pC_{J,J}^{i,j})_{1\leq i,j\leq k})$ for $J=1,2$
(resp. $\det((\pC_{J,J}^{i,j})_{2\leq
i,j\leq k}$ for $J\geq 3$)
by the same permutation
$\sigma\in S_k$ in the natural way.
\end{proposition}
\begin{proof}
It follows from Lemmas \ref{L:locate-min} and
\ref{L:minimal-weight} above.
\end{proof}

\subsection{Local at each pole}

For ease of notation, we drop the
subindex $J$ for the rest of this subsection.
One should understand that
$d,A_i,F_{ip-j},D_n$ stand for $d_J, A_{J,i}, F_{J,ip-j},
D_{J}^{[n]}$,
respectively.
Let $1\leq n\leq d-1$ and let $S_n$ be the permutation group.
Let $D_n := \det
(F_{ip-j})_{1\leq i,j \leq n} \in\cO_1[A_1,\ldots,A_d]$.
Then we have the formal expansion of $D^{[n]}$:
\begin{eqnarray*}
D^{[n]}=\sum_{\sigma\in S_n} \sgn(\sigma) \sum\prod_{i=1}^{n}
g_{\sigma, i},
\end{eqnarray*}
where the second $\sum$ runs over all
terms $g_{\sigma, i}$ of the polynomial
$F_{ip-\sigma(i)}$ in $\cO_1[A_1,\ldots,A_d]$.

\begin{proposition}\label{P:keyproposition}
Let $1\leq n\leq d$.
Then there is a unique monomial in the above formal expansion of
$D^{[n]}$ with highest lexicographic order (according to
$A_d,\ldots,A_1$). Moreover, the $p$-adic order of this monomial
(with coefficient)
is minimal among the $p$-adic orders of all monomials in the above
formal expansion.
\end{proposition}

\begin{remark}\label{R:mini-valuation}
We shall fix the unique $\sigma_0$ found in
the proposition for the rest of the paper.
The minimal $p$-adic order of this monomial
(with coefficient) is equal to (\ref{E:s_J})
while every row achieve its minimal order
in Lemma \ref{L:F&H}c).
We shall use this fact later.
\end{remark}

\begin{proof}
Denote by $r$ the least non-negative residue of $p \mod d$. Recall
the $n$ by $n$ matrix $\br_n:=\{r_{ij}\}_{1\leq i,j\leq n}$ where
$r_{ij}:= d\lceil\frac{ri-j}{d}\rceil- (ri-j)$. The properties of
this matrix can be found in \cite[Lemma 3.1]{Zhu:1}. Let
$\prod^n_{i=1}h_{\sigma, i}$ be a
highest-lexicographic-order-monomial in the formal expansion of
$D^{[n]}$. Then $h_{\sigma, i}$ must be the
highest-lexicographic-order-monomial in $F_{ip-\sigma(i)}$,
which is easily seen to be
$c_{ip-\sigma(i)}A_d^{\frac{ip-\sigma(i)}{d}}$ or
$c_{ip-\sigma(i)}A_d^{\lfloor\frac{ip-\sigma(i)}{d}\rfloor}A_{d-r_{i,
\sigma(i)}}$ depending on $r_{i, \sigma(i)}=0$ or not, where
$c_{ip-\sigma(i)}\in \Omega_1$. We show first that
\begin{eqnarray} \label{perm 0:eq}
\sigma(i)=k \mbox{ for any } 1\le i, k\le n \mbox{ with } r_{i
k}=0.
\end{eqnarray}

Suppose that (\ref{perm 0:eq}) does not hold. Pick a pair $(i,k)$
among the pairs failing (\ref{perm 0:eq}) such that
$|\sigma(i)-k|$ is minimal. Say $\sigma(j)=k$. Define another
permutation $\sigma'\in S_n$ by $\sigma'(i)=k$ and
$\sigma'(j)=\sigma(i)$ while $\sigma'(s)=\sigma(s)$ for all other
$s$. Denote by $h_{\sigma', t}$ the
highest-lexicographic-order-monomial in $F_{tp-\sigma'(t)}$. Then
it is easy to see that the lexicographic order of
$\prod^n_{t=1}h_{\sigma', t}$ is strictly higher than that of
$\prod^n_{t=1}h_{\sigma, t}$, which is a contradiction. Therefore
(\ref{perm 0:eq}) holds.

Notice that for permutations $\sigma''\in S_n$ satisfying
(\ref{perm 0:eq}), the degree of $A_d$ in
$\prod^n_{t=1}h_{\sigma'', t}$ does not depend on the choice
of $\sigma''$, where $h_{\sigma'', t}$ is the
highest-lexicographic-order-monomial in $F_{tp-\sigma''(t)}$. Then
the proof of \cite[Lemma 3.2]{Zhu:1} shows that there exists a
unique $\sigma_0\in S_n$ such that $\prod^n_{t=1}h_{\sigma_0, t}$
has highest lexicographic order among the corresponding monomials
for all $\sigma''\in S_n$ satisfying (\ref{perm 0:eq}). In fact,
$\sigma_0$ is exactly the permutation in \cite[Lemma 3.2]{Zhu:1}.
By the above discussion, this monomial $\prod^n_{t=1}h_{\sigma_0,
t}$ also has the unique highest lexicographic order in the formal
expansion of $D^{[n]}$.

Next we show that
 the $p$-adic order of $\prod^n_{i=1}h_{\sigma_0, i}$ is minimal
 (among the $p$-adic orders of the monomials in the formal
expansion of $D^{[n]}$). Let $\prod^n_{i=1}g_{\sigma, i}$ be an
arbitrary monomial in the formal expansion of $D^{[n]}$. Then
clearly $\ord_p(g_{\sigma, i})\ge
\lceil\frac{ip-\sigma(i)}{d}\rceil=\frac{pi-\sigma(i)+r_{i,
\sigma(i)}}{d}$ for all $1\le i\le n$. Since $h_{\sigma_0, i}$ is
the highest-lexicographic-order-monomial
 in $F_{ip-\sigma_0(i)}$, one sees easily that $\ord_p(h_{\sigma_0,
i})=\frac{pi-\sigma_0(i)+r_{i, \sigma_0(i)}}{d}$ for all $1\le
i\le n$. From (\ref{perm 0:eq}) it is easy to see that $r_{i, j}-r_{i,
\sigma_0(i)}\ge j-\sigma_0(i)$ for all $1\le i, j\le n$. It
follows that $\ord_p(\prod^n_{i=1}g_{\sigma, i})\ge
\ord_p(\prod^n_{i=1}h_{\sigma_0, i})$.
\end{proof}

\begin{remark}\label{R:f_n}
In the proof of Proposition \ref{P:keyproposition}
we have noticed that the $\sigma_0$ is exactly
the permutation in \cite[Lemma 3.2]{Zhu:1}.
Therefore,
one can always take $t_0=0$, that is,
$f_n^0(\vec{A})\neq 0$ in \cite[Lemma 3.5]{Zhu:1}.
\end{remark}

\section{Newton polygon of $\alpha_1$}

Recall that $\HP(\Ae)$ lives on the real plane over the interval
$[0,d]$.
Because of Remark \ref{R:1}, one only has to consider the
part of $\NP_p(f)$ with slope $<1$, that is, to consider
the part of
$\NP_p(f)$ over the interval $[0,d-\ell]$. This part is our focus of
this section.
Suppose for some $1\leq k\leq d-\ell$, the point
$(k,c_0)$ is a vertex on $\HP(\Ae)$.
Then one notices that
$$
c_0=\sum_{J=1}^{2}\sum_{i=1}^{k_J}i/d_J +\sum_{J=3}^{\ell}\sum_{i=1}
^{k_J}(i-1)/d_J
$$
for a sequence of nonnegative integers
$k_1,\ldots,k_\ell$
such that $k_1+\ldots+k_\ell=k$.
This sequence is unique because $(k,c_0)$ is a vertex.
{}From now on we fix such a $k$.

For our purpose we also fix the residue classes of $p\bmod d_J$
for all $J$. Let $r_{J,ij}$ be the least nonnegative
residue of $-(ip-j)\bmod d_J$.
Let $\sigma_0$ be the permutation in $S_k$ which is
the union of those permutations found in Proposition
\ref{P:keyproposition} locally at each pole $P_J$.
Let $s_{J_1}$ be the rational number defined by
\begin{eqnarray}\label{E:s_J}
s_{J_1} &:=&
\frac{(p-1)k_{J_1}(k_{J_1}\pm 1)/2}{d_{J_1}} +
\frac{\sum_{i=1}^{k_{J_1}}r_{J_1, i,\sigma_0(i)}}{d_{J_1}}
\end{eqnarray}
where $+$ and $-$ is taken according to
$J_1=1,2$ or $J_1\geq 3$.
Let $s_0:=s_1+\cdots +s_\ell$.
Clearly $s_0-c_0(p-1)<k\leq d-\ell$.

Let $\alpha$ and $\beta$ be
the integers chosen in Lemmas \ref{L:compare}
and \ref{L:finite-term}
(they depend only on $d$).
Let $\Qe':=\Qe(\gamma^{1/d_1},\ldots, \gamma^{1/d_\ell})$.

\begin{lemma}\label{L:w}
For any $J_1=1,2$, $1\leq J\leq \ell$ and for any $n,i$ in their
range, for $p$ large enough, there is a polynomial
$G_{J_1,J}^{n,i}$ in $\Qe'(\vec{P})[\vec{A}]$ such that
$$
{}^\alpha K_{J_1,J}^{n,i}
= \gamma^{(p-1)n/d_{J_1}} U_{J_1,n} G_{J_1,J}^{n,i}
\bmod \gamma^{(p-1)n/d_{J_1}+d+1}.
$$
For the case $J_1\geq 3$,
one has a similar $G_{J_1,J}^{n,i}$
such that
$$
{}^\alpha K_{J_1,J}^{n,i}
= \gamma^{(p-1)(n-1)/d_{J_1}} U_{J_1,n} G_{J_1,J}^{n,i}
\bmod \gamma^{(p-1)(n-1)/d_{J_1}+d+1},
$$
where $U_{J_1,n}$ is a $p$-adic unit depending only on
the the row index $(J_1,n)$.
\end{lemma}
\begin{proof}
We use the same technique as \cite{Zhu:1}, so we only outline our
proof here for the case $J_1=J=1$. Let $n_j\in [-\alpha,\alpha]$
for $j\neq J_1$. Let $n_{J_1} =np+\sum_{j\neq J_1}n_j - i$. For
any $\vec{n} =(n_1,\ldots,n_\ell)$ in this range, we have
$$
F_{j,n_j} \equiv \gamma^{\frac{n_j}{d_j}}Q_j \bmod \gamma^{
\frac{n_j}{d_j}+d+1}
$$
and
$$
F_{J_1,n_{J_1}}
 =\gamma^{\frac{n_{J_1}}{d_{J_1}}}
  V_{J_1} Q_{J_1} \bmod \gamma^{\frac{n_{J_1}}{d_{J_1}}+d+1},
$$
where $Q_j$'s and $Q_{J_1}$ are in $\Qe'[\vec{A}]$ independent of
$p$ and $V_{J_1}$ is some $p$-adic unit depending only on
the row index $J_1$. Now let $\vec{n}$ be in the range for ${}^\alpha
K_{J_1,J}^{n,i}$ such that $n_j$'s vary in $[-\alpha,\alpha]$ and
$$
\frac{i}{d_J}-\frac{n}{d_{J_1}} +\frac{np+\sum_{j\neq J_1}n_j -
i}{d_{J_1}} + \sum_{j\neq J_1}\frac{n_j}{d_j} \geq
(p-1)n/d_{J_1}.
$$
Then by the formula of Lemma \ref{L:Hjj} (1),  and for
$p$ large enough,
$$
{}^{\alpha} K_{J_1,J}^{n,i}
= \gamma^{\frac{i}{d_J}-\frac{n}{d_{J_1}}}
  {}^\alpha H_{J_1,J}^{np,i}
\equiv
  \gamma^{\frac{n(p-1)}{d_{J_1}}} W G_{J_1,J}^{n,i}
\bmod \gamma^{(p-1)n/d_{J_1}+d+1},
$$
where $W$ is a suitable $p$-adic unit.
The rest of the cases are similar.
\end{proof}

\begin{proposition}\label{P:C}
Let notation be as in Lemma \ref{L:w}.
Let $\bK   :=  ({}^{\alpha}K_{J_1,J}^{n,i})$.
For $p$ large enough,
there are a polynomial $Y_k$ in $\Qe'(\vec{P})[\vec{A}]$
and some $p$-adic unit $U$
such that
$$
\det \bK^{[k]}
\equiv
\gamma^{c_0(p-1)}U
Y_k \bmod \gamma^{c_0(p-1)+d+1}.
$$
\end{proposition}
\begin{proof}
By Lemmas~\ref{L:w} and \ref{L:F&H}(c) we have
\begin{eqnarray*}
\det \bK^{[k]}\equiv \gamma^{c_0(p-1)}U
\det  \bG^{[k]} \mod \gamma^{c_0(p-1)+d+1},
\end{eqnarray*}
where $\bG^{[k]}$ is the matrix we obtain via replacing
${}^{\alpha}K_{J_1,J}^{n,i}$ by $G_{J_1,J}^{n,i}$ in $\bK^{[k]}$,
and $U$ is the product of the $U_{J_1, n}$'s for the pairs $(J_1,
n)$ whose corresponding row appears in $\bM^{[k]}$. Now just set
$Y_k=\det \bG^{[k]}$.
\end{proof}

\begin{lemma}\label{L:keylemma}
Let $1\leq k \leq d-\ell$.
(1) For $p$ large enough one has
\begin{eqnarray}\label{E:compare2}
\ord_p (\det\bM^{[k]} - \det{\pM}^{[k]}) & > & \frac{s_0}{p-1}.
\end{eqnarray}
(2)
Let $\alpha$ and $\beta$ be
the integers chosen in Lemmas \ref{L:compare}
and \ref{L:finite-term}
(they depend only on $d$).
Then
\begin{eqnarray}
\ord_p (\det\pM^{[k]} - \det \bK^{[k]})
&> & \frac{s_0}{p-1}.
\end{eqnarray}
\end{lemma}

\begin{proof}
(1) Note that $d\geq k > s_0 - c_0(p-1)$. Note that in $\pM^{[k]}$
the row minimal $p$-adic order is the same as that for
$C_{J_1,J}^{n,i}$ in Lemma \ref{L:F&H} (c). By Lemma
\ref{L:compare}, for $p$ large enough one has
\begin{eqnarray}
\ord_p(\det \bM^{[k]} - \det \pM^{[k]})
 & > &
c_0 + \frac{d}{p-1}\geq \frac{s_0}{p-1}.
\end{eqnarray}

(2)
By Lemma \ref{L:finite-term}, one knows that
$\ord_p({C'}_{J_1,J}^{n,i} - K_{J_1,J}^{n,i})
> \frac{n-1}{d_{J_1}}+ \frac{d}{p-1}$.
Thus
\begin{eqnarray}
\ord_p(\det\pM^{[k]} - \det \bK^{[k]})
 & > &
c_0 + \frac{d}{p-1} \geq \frac{s_0}{p-1},
\end{eqnarray}
since $s_0 -c_0(p-1)<k$.
\end{proof}

In any formal expansion we group the terms with same $p$-adic
orders together and then write this in increasing order. For any
number $t$ in $\Qe$ if a term can be written as $\gamma^{t}u$ for
some $u$ with $\ord_pu=0$, then $u$ is called the
$\gamma^t$-coefficient of this term. Let $\bM(\hf)$ denote the
specialization of $\bM$ at variables $\vec{A}$ by assigning
$\vec{A}$ as the Teichm\"uller lifts of coefficients of $f\bmod
\cP$ (see \cite[Section 1]{Zhu:2} for more details).

\begin{proposition}\label{P:F_p}
Let $1\leq k\leq d-\ell$. Let  $(k, c_0)\in \Re^2$ be a vertex of
the slope $<1$ part of $\HP(\Ae)$, where $1\le k\le d-\ell$. There
is a Zariski dense open subset $\cU_k$ defined over $\Qe$ in $\Ae$
such that if $f\in\cU_k(\bar\Qe)$ and if $\cP$ is a prime ideal in
the ring of integers of $\Qe (f)$ lying over $p$, one has
$\lim_{p\rightarrow\infty}\ord_p\det(\bM(\hf)^{[k]})=c_0$.
\end{proposition}

\begin{proof} Without loss of generality, we fix the residues of
$p$ as above. Consider the $\gamma$-expansions of $\det\bM^{[k]},
\det\pM^{[k]}$, and $\det{\bK}^{[k]}$. By Lemma \ref{L:keylemma}, their
$\gamma^{s_0}$-coefficients are the same. Proposition \ref{P:C}
implies that for $p$ large enough there is a polynomial $G$ in
$\Qe(\vec{P})[\vec{A}]$ such that the $\gamma^{s_0}$-coefficient
is congruent to $U\,G\bmod\gamma$ for some $p$-adic unit $U$.
Moreover, from the proofs of Lemma \ref{L:w} and Proposition
\ref{P:C}, one observes easily that the monomials of $G$ are a
subset of all monomials in the formal expansion of $\det\pM^{[k]}$
(with all $\gamma^\Qe$-factors squeezed out from its coefficients
at appropriate places).

We claim that the $\gamma^{s_0}$-coefficient in $\det \pM^{[k]}$
is nonzero because it has a unique monomial (in variable
$\vec{A}$) among all monomials of minimal weight in its formal
expansion. We first look locally at an arbitrary pole $P_J$ where
$1\leq J\leq \ell$. By Proposition \ref{P:keyproposition} there is
a unique local monomial among all terms in $\det D_J^{[k_J]}$
for $J=1,2$ and $\det D_J^{[k_J-1]}$ for $J\geq 3$.
This local monomial corresponds to a permutation $\sigma_{J,0}\in
S_{k_J}$.
Note that the composition of these $\sigma_{J,k_J}$'s  for all $J$ is
equal to $\sigma_0$ defined in the beginning of the section. Then
the unique monomial we desire is precisely the product of these
local monomials (see Lemma \ref{L:minimal-weight} and Proposition
\ref{P:weight}).
 By the remark in last paragraph, it is not hard to see
that $G\neq 0$.

Let $\gamma^{>s_0}$ denote all those terms with $p$-adic order
$>\frac{s_0}{p-1}$. Recall from Lemma \ref{L:keylemma} and
Proposition \ref{P:C} that one has the $p$-adic unit
$U$ (as in the above paragraphs) and some polynomials
$G'_m$ and $G'$ (in $\Qp(\vec{\hP})[\vec{A}]$) such that
$$\det(\bM^{[k]}) =
\sum_{c_0\leq m < s_0}\gamma^{m}U\, G'_m
+\gamma^{s_0} U\, G' +\gamma^{>s_0}
$$
and $G'\equiv G\bmod \gamma$ for the polynomial $G$ (same
$G$ as in above paragraphs) in $\Qe(\vec{P})[\vec{A}]$ independent
of $p$. If $G(f)\not\equiv 0\bmod \cP$ (the specialization of $G$
at $f$ over $\Qe(\vec{P})$) then $\ord_p(G'(\hf))=0$. For $m<s_0$
one has $\ord_pG_m'(\hf)=0$ or $\geq 1$. Thus if $G(f)\neq 0$ then
for $p$ large enough one has $ c_0\leq \ord_p(\det\bM(\hf)^{[k]})
\leq \frac{s_0}{p-1}.$ But we already know from the beginning of
this section that $0\leq \frac{s_0}{p-1}-c_0\leq
\frac{d-\ell}{p-1}$ and hence by simple calculus one has that
$\lim_{p\rightarrow \infty} \ord_p(\det\bM(\hf)^{[k]}) = c_0$.

Last, taking the norm of $G$ from $\Qe(\vec{P})[\vec{A}]$ to
$\Qe[\vec{A}]$ with the automorphism acting on $\vec{A}$
trivially, one gets a polynomial $g$ in $\Qe[\vec{A}]$. Let $\cV$
be the complement of the variety defined by $g=0$ in $\Ae$. It
is Zariski dense in $\Ae$ because $g\neq 0$.
\end{proof}

\section{A transformation theorem
from Newton polygons of $\alpha_1$ to $\alpha_a$}

We refer the reader to \cite{Serre62, Monsky70} for basic facts
about Serre's theory of completely continuous maps and Fredholm
determinants.  Let $\Cp$ be the $p$-adic completion of $\bar\Qp$.
For any $\Ce_p$-Banach spaces $E$ and $F$ that admit orthonormal
bases, denote by $\mathscr{C}(E, F)$ the set of completely
continuous $\Ce_p$-linear maps from $E$ to $F$. We say that a
matrix $M$ over $\Ce_p$ is {\it nuclear} if there exist a
$\Ce_p$-Banach space $E$ and a $u\in \mathscr{C}(E,\,E)$ such that
$M$ is the matrix of $u$ with respect to some orthonormal basis of
$E$.
If $M=(m_{ij})_{i,\, j\ge 1}$ is a matrix over $\Cp$, then $M$ is
nuclear if and only if $\lim_{i\to \infty}(\inf_{j\ge
1}\ord_pm_{i, \, j})=+\infty$. Recall $\ord_q (\cdot) =
\ord_p(\cdot)/a$ for $q=p^a$.

\begin{lemma} \label{Ma:lemma}
Let
$\vec{M}=(M_0,\, M_1, \cdots,
M_{a-1})$ be an $a-$tuple of nuclear matrices over $\Ce_p$. Set
\begin{eqnarray*}
\vec{M}_{[a]}:=
\begin{pmatrix}
0      &        &\cdots  & 0       & M_{a-1} \\
M_0    &  0     &        &         & 0       \\
0      & M_1    & 0      &         & \vdots  \\
\vdots &        & \ddots & 0       &         \\
0      & \cdots & 0      & M_{a-2} & 0
\end{pmatrix}.
\end{eqnarray*}
Then $\det(1-(M_{a-1}\cdots M_1M_0)T^a)=\det(1-\vec{M}_{[a]}T)$.
\end{lemma}

Lemma~\ref{Ma:lemma} follows directly from

\begin{lemma} \label{ua:lemma}
Let $\{E_i\}_{i\in \Ze/a\Ze}$ be a family of Banach spaces over
$\Ce_p$ that admit orthonormal basis. Set $E=E_0\oplus E_1\oplus
\cdots \oplus E_{a-1}$ equipped with the supremum norm, that is for
$v=(v_0,\ldots,v_{a-1})$ in $E$ one has $||v||=\max_{i=0}^{a-1}||v_i||$,
where $||\cdot||$ are the norms on $E$ and $E_i$'s, respectively.  Let
$u_i\in \mathscr{C}(E_i, \, E_{i+1})$ and set $u\in \mathscr{C}(E, \,
E)$ such that $u|_{E_i}=u_i$. Then
\begin{eqnarray*}
\det(1-(u_{a-1}\cdots u_1u_0)T^a)=\det (1-u T).
\end{eqnarray*}
\end{lemma}
\begin{proof}
By \cite[page 77,\, Corollaire 3]{Serre62} we have $\det(1-u
T)=\exp (-\sum^{\infty}_{s=1}\Tr(u^s)T^s/s)$. Notice that for any
$s\in \Ze_{\ge 1}$, the trace $\Tr((u_{i+a-1}\cdots
u_{i+1}u_i)^s)$ is independent of $i\in \Ze/a\Ze$. Clearly
$\Tr(u^s)=0$ unless $a|s$. Thus
\begin{eqnarray*}
\det(1-u T)&=&\exp (-\sum^{\infty}_{s=1}\Tr(u^s)T^s/s) =\exp
(-\sum^{\infty}_{s=1}\Tr(u^{as})T^{as}/(as)) \\
&=&\exp (-\sum^{\infty}_{s=1}\sum_{i\in
\Ze/a\Ze}\Tr((u_{i+a-1}\cdots u_{i+1}u_i)^s))T^{as}/(as)) \\
&=&\exp (-\sum^{\infty}_{s=1}\Tr((u_{a-1}\cdots
u_{1}u_0)^s))T^{as}/s) \\
&=&\det(1-(u_{a-1}\cdots u_1u_0)T^a).
\end{eqnarray*}
This concludes our proof.
\end{proof}

\begin{remark} \label{Ma:remark}
 Lemmas~\ref{Ma:lemma} and \ref{ua:lemma} still hold when $\Ce_p$ is
replaced by any field $K$ equipped with a nontrivial complete
non-Archimedean valuation.  But we shall not need this more general
fact in the present paper.
\end{remark}

For any nuclear matrix $M=(m_{ij})_{i, j\ge 1}$ and $k\in \Ze_{\ge
1}$, denote by $M^{[k]}$ the $k\times k$ submatrix of $M$ consisting of its first
$k$ rows and columns.

\begin{proposition} \label{NP Ma:prop}
Let $M=(m_{ij})_{i, j\ge 1}$ be a nuclear matrix over $\Ce_p$
and let $g\in \Gal(\bar{\Qe}_p/\Qe_p)$. Fix $k\in \Ze_{\ge 1}$ and
denote by $C_k$ the coefficient of $T^k$ in
$\det(1-(M^{g^{a-1}}\cdots M^gM)T)$. Denote by $\mathscr{A}$ the set of
$k\times k$ submatrices of $M$ contained in the first $k$ rows of
$M$, and denote by $\mathscr{B}$ the set of all other $k\times k$ submatrices
of $M$.
Set $t_{\mathscr{A}}=\inf_{W\in \mathscr{A}} \ord_p\det W$ and
$t_{\mathscr{B}}=\inf_{W\in \mathscr{B}}\ord_p\det W$.
Consider the following conditions:
\begin{itemize}
\item[(i)] $2\ord_p\det M^{[k]}<t_{\mathscr{A}}+t_{\mathscr{B}}$;
\item[(ii)] $2\ord_qC_k<t_{\mathscr{A}}+t_{\mathscr{B}}$ and
$t_{\mathscr{A}}<t_{\mathscr{B}}$;
\item[(iii)] $\ord_qC_k=\ord_p\det M^{[k]}$.
\end{itemize}
Then (i)$\iff$(ii)$\Longrightarrow$(iii).
\end{proposition}
\begin{proof}
Notice that $\ord_p\det M^{[k]}\ge t_{\mathscr{A}}$.
So (i) is equivalent to
\begin{eqnarray} \label{A, B, M:equa}
\min (\frac{t_{\mathscr{A}}+t_{\mathscr{B}}}{2},\,
t_{\mathscr{B}})>\ord_p\det M^{[k]}.
\end{eqnarray}
It suffices to show that (ii)$\Rightarrow$(\ref{A, B,
M:equa})$\Rightarrow$(iii).  Let $\vec{M}:=(M, \, M^g, \, \cdots,
\,M^{g^{a-1}})$. Then we have $\vec{M}_{[a]}$ in Lemma~\ref{Ma:lemma}
and $\det(1-\vec{M}_{[a]}T)=\det(1-(M^{g^{a-1}}\cdots M^gM)T^a)$. Thus
$C_k$ is the coefficient of $T^{ak}$ in $\det(1-\vec{M}_{[a]}T)$,
which is the infinite sum of $(-1)^{ak}\det N$ for $N$ running over
all principal $ak\times ak$ submatrices of $\vec{M}_{[a]}$. Let $N$ be
such a matrix, and let $N_s$ be the intersection of $N$ and $M^{g^s}$
as submatrices of $\vec{M}_{[a]}$ for all $0\le s\le a-1$. It is easy
to see that $\det N=(-1)^{(a-1)k}\prod_{0\le s\le a-1}\det N_s$ or $0$
depending on whether every $N_s$ is a $k\times k$ matrix or not.  So
we may assume that every $N_s$ is a $k\times k$ matrix.  Think of
$N_s$ as a submatrix of $M^{g^s}$ from now on.  Let $X=\{s: 0\le s\le
a-1 \mbox{ and } (N_s)^{g^{-s}}\in \mathscr{A}\setminus \{M^{[k]}\}\}$
and $Y=\{s: 0\le s\le a-1 \mbox{ and } (N_s)^{g^{-s}}\in
\mathscr{B}\}$.  We shall think of the families $\{M^{g^s}\}_{0\le
s\le a-1}$ and $\{N_s\}_{0\le s\le a-1}$ as parameterized by
$\Ze/a\Ze$.  Then $X$ and $Y$ are disjoint subsets of $\Ze/a\Ze$.
Since $N$ is principal, the set of the columns of $N_s$ as a subset in
$\Ze_{\ge 1}$ is exactly the same as the set of the rows of $N_{s-1}$.
Consequently, if $s\in X$, then $s-1\in Y$.  Let $Y'=\{s-1:s\in X\}$
and $Z=(\Ze/a\Ze)\setminus (X\cup Y)$.  Then $\Ze/a\Ze$ is the
disjoint union of $X\cup Y', \, Y\setminus Y'$ and $Z$.  If $s\in X$,
then $\ord_p(\det N_s\cdot \det N_{s-1})\ge
t_{\mathscr{A}}+t_{\mathscr{B}}$.  If $s\in Y\setminus Y'$, then
$\ord_p \det N_s\ge t_{\mathscr{B}}$.  If $s\in Z$, then $\ord_p\det
N_s= \ord_p\det M^{[k]}$.  Therefore
\begin{eqnarray} \label{det N:equa}
\ord_q\det N\ge \min (\frac{t_{\mathscr{A}}+t_{\mathscr{B}}}{2},\,
t_{\mathscr{B}},\,\ord_p\det M^{[k]}),
\end{eqnarray}
and hence
\begin{eqnarray} \label{Ck:equa}
\ord_qC_k\ge \min (\frac{t_{\mathscr{A}}+t_{\mathscr{B}}}{2},\,
t_{\mathscr{B}},\,\ord_p\det M^{[k]}).
\end{eqnarray}

(\ref{A, B, M:equa})$\Rightarrow$(iii):
Clearly there is a unique $N$ with $X=Y=\emptyset$, i.e.
$(N_s)^{g^{-s}}=M^{[k]}$ for all $0\le s\le a-1$.
Denote it by $\cN$.
We have $\ord_q\det \cN=\ord_p\det M^{[k]}$.
If $N\neq \cN$, then $X$ or $Y\setminus Y'$ is nonempty and hence
from (\ref{A, B, M:equa}) and the derivation of (\ref{det N:equa})
we see that $\ord_q\det N>\ord_p\det M^{[k]}$. Now (iii) follows immediately.

(ii)$\Rightarrow$(\ref{A, B, M:equa}): (\ref{A, B, M:equa}) follows
directly from (ii) and (\ref{Ck:equa}).
\end{proof}

\begin{theorem} \label{T:F_q}
Let $M, \, g, \, k\,$ and $C_k$ be as in Proposition~\ref{NP
Ma:prop}. Let $h_1\le h_2\le \cdots$ be a non-decreasing sequence
in $\Re$ satisfying $h_i\le \inf_{j\ge 1}\ord_pm_{ij}$ for all
$i\ge 1$. Consider the following conditions:
\begin{itemize}
\item[(i)] $\ord_p\det
M^{[k]}<\sum_{1\le i\le k}h_i+\frac{h_{k+1}-h_k}{2}$;
\item[(ii)] $\ord_qC_k<\sum_{1\le i\le
k}h_i+\frac{h_{k+1}-h_k}{2}$;
item[(iii)] $\ord_qC_k=\ord_p\det
M^{[k]}$.
\end{itemize}
Then (i)$\iff$(ii)$\Longrightarrow$(iii).
\end{theorem}
\begin{proof}
Let $t_{\mathscr{A}}$ and $t_{\mathscr{B}}$ be as in
Proposition~\ref{NP Ma:prop}.  Then $\sum_{1\le i\le
k}h_i+\frac{h_{k+1}-h_k}{2}\le
\min(\frac{t_{\mathscr{A}}+t_{\mathscr{B}}}{2},\, t_{\mathscr{B}})$.
So (i) follows from (ii) and (\ref{Ck:equa}).  Thus
Theorem~\ref{T:F_q} follows from Proposition~\ref{NP Ma:prop}.
\end{proof}

\begin{remark}
Theorem \ref{T:F_q} is a Wan-type theorem in relating the Newton
polygon to its tight lower bound Hodge polygon: In \cite[Theorem
8]{Wan:1}, Wan showed that the Newton polygon for $\alpha_1$ (more
precisely, the Fredholm determinant of the nuclear matrix
representing $\alpha_1$ with respect to the specific basis)
coincides with the Hodge one if and only if the Newton polygon for
$\alpha_a$ does. Our result in Theorem~\ref{T:F_q} generalizes it
and says that the Newton polygon for $\alpha_1$ is close to the
Hodge one if and only if the Newton polygon for $\alpha_a$ is.
\end{remark}

\begin{proof}[Proof of Theorem \ref{T:main}]
For any vertex $(k, c_0)\in \Re^2$ (but not the right end point)
of the slope $<1$ part of $\HP(\Ae)$, where $1\le k\le d-\ell$,
let $\cU_k$ be the Zariski dense open subset in
Proposition~\ref{P:F_p}. Let $f\in\cU_k(\bar\Qe)$. Then
$\lim_{p\rightarrow\infty}\ord_p\det(\M(\hf)^{[k]})=c_0$. Recall
$\phi(\cdot)$ from the beginning of section~\ref{S:3.1}. Say the
coefficients of $f \mod \cP$ lie in $\mathbb{F}_{p^a}$.  Set
$M:=\M(\hf)$ and $h_i:=\phi(e_i)$ for all $i\ge 1$ in
Theorem~\ref{T:F_q}. Notice that $\sum_{1\le i\le k}h_i=c_0$.
Since $(k, c_0)$ is a vertex of $\HP(\Ae)$, we have $h_{k+1}>
h_k$. In particular, when $p$ is large enough, we have $\ord_p\det
M^{[k]}< c_0 + \frac{h_{k+1}-h_k}{2}$. Combining this with
Lemma~\ref{L:F&H}(c), one observes that the hypotheses of
Theorem~\ref{T:F_q} are satisfied. Recall the maps $\alpha_1$ and
$\alpha_a$ defined in Lemma 2.9 and section 2.5 of \cite{Zhu:26}.
These maps are not the same as the maps defined in section 2 of
this article, but are the specialization of those maps in section
2 at the Teichm\"uller lifts of coefficients of $f\bmod \cP$. Then
$M^{\tau^{-1}}$ and $M^{\tau^{-1}}\cdots
M^{\tau^{-(a-1)}}M^{\tau^{-a}}$ are the matrices of $\alpha_1$ and
$\alpha_a$ (over $\Omega'_a$) with respect to the formal basis
$\vec{b}_{\rm w}=\{1, Z^i_1, \cdots, Z^i_{\ell}\}_{i\ge 1}$ of
$\cH$ respectively. Notice that $M^{\tau^a}=M$. By
Theorem~\ref{T:F_q} one has $\lim_{p\rightarrow\infty}\ord_qC_k =
c_0$, where $C_k$ is the coefficient of $T^k$ in
$\det(1-(M^{\tau^{a-1}}\cdots
M^{\tau}M)T)=\det(1-(M^{\tau^{-1}}\cdots
M^{\tau^{-(a-1)}}M^{\tau^{-a}})T)=\det_{\Omega'_a}(1-\alpha_aT)$.
Set $\cU$ to be the intersection of $\cU_k$ for all such vertices
$(k, c_0)$. Then for any $f\in\cU(\bar\Qe)$, we have
$\lim_{p\rightarrow\infty} \NP_q(\det_{\Omega'_a}(1-\alpha_a
T)\mod T^{d-\ell+1})=\HP(\Ae)$. Now Theorem \ref{T:main} follows
from Remark~\ref{R:1} and the fact that the slope $<1$ part of
$\NP_p(f)$ coincides with $\NP_q(\det_{\Omega'_a}(1-\alpha_a
T)\mod T^{d-\ell+1})$ (see \cite[Proposition 2.10]{Zhu:26}).
\end{proof}

\begin{remark}
(1) Our main result Theorem \ref{T:main} is related but not
included in a conjecture of Daqing Wan (see \cite[Conjectures 1.12
and 1.14]{Wan:2}).

(2) This paper is concerned with the space of all
one-variable rational function
with fixed poles on the projective line. One naturally wonders if
there is a multivariable generalization of Theorem \ref{T:main}.
We do not know the answer.

\end{remark}


\begin{thebibliography}{99}

\bibitem{Berthelot}
{\sc Pierre Berthelot:} Cohomologie rigide et th\'eorie de Dwork:
le ca des sommes exponentielles, in {\it Cohomologie p-adique},
Soci\'et\'e Math\'ematique de France, Ast\'erisque {\bf 119--120}
(1984), 17-49.

\bibitem{BGR:1}
{\sc S. Bosch; U. Guntzer; R. Remment:} Non-Archimedean analysis,
{\it Grundlehren der Mathematischen Wissenschaften} Vol. {\bf
261}, Springer-Verlag, Berlin, 1984.

\bibitem{Crew:1}
{\sc Richard Crew:}
Etale $p$-covers in characteristic $p$,
{\it Compositio Math.},
{\bf 52} (1984),
31--45.

\bibitem{Katz:1}
{\sc Nicholas M. Katz:}
Gauss sums, Kloosterman sums, and monodromy groups,
{\it Annals of mathematics studies} vol {\bf 116},
Princeton University Press, 1988.

\bibitem{Monsky70}
{\sc Paul Monsky:} $p$-adic analysis and zeta functions, Lectures
in Mathematics, Department of Mathematics, Kyoto University,
Kinokuniya Book-Store Co., Ltd., Tokyo, 1970.

\bibitem{Robba:1}
{\sc Philippe Robba:}
Index of $p$-adic differential operators III.
Application to twisted exponential sums, in {\it Cohomologie
p-adique},
Soci\'et\'e Math\'ematique de France, Ast\'erisque
{\bf 119--120} (1984),
191--266.

\bibitem{SZ:1}
{\sc Jasper Scholten; Hui June Zhu:}
Hyperelliptic curves in characteristic $2$,
{\it Math. Research Letters}
{\bf 17} (2002), 905--917.

\bibitem{Serre62}
{\sc Jean-Pierre Serre:}
Endomorphismes compl\`etement continus des espaces de Banach $p$-adiques,
{\it Inst. Hautes \'Etudes Sci. Publ. Math.} {\bf 12} (1962), 69--85.

\bibitem{Wan:1}
{\sc Daqing Wan:}
Newton polygons of zeta functions and $L$-functions
{\it Ann. Math.}
{\bf 137}
(1993), 247--293.

\bibitem{Wan:2}
{\sc Daqing Wan:} Variation of Newton polygons for $L$-functions
of exponential sums.
{\it Asian J. Math.} {\bf 8}  (2004), 427--474.

\bibitem{Wan:3}
{\sc Daqing Wan:}
Rank one case of Dwork's conjecture.
{\it J. of Amer. Math. Soc.}
{\bf 13} (2000), 853--908.

\bibitem{Wan:4}
{\sc Daqing Wan:}
Higher rank case of Dwork's conjecture.
{\it J. of Amer. Math. Soc.}
{\bf 13} (2000), 807--852.

\bibitem{Zhu:1}
{\sc Hui June Zhu:} $p$-adic variation of $L$ functions of one
variable exponential sums, I. {\it Amer. J. Math.}
{\bf 125} (2003), 669--690.

\bibitem{Zhu:2}
{\sc Hui June Zhu:}
Asymptotic variation of $L$ functions of
one-variable exponential sums.
{\it J. Reine Angew. Math.}
{\bf 572} (2004), 219--233.

\bibitem{Zhu:26}
{\sc Hui June Zhu:} $L$ functions of exponential sums over one
dimensional affinoids: Newton over Hodge. {\it Inter. Math. Res.
Notices.}, vol 2004, no. 30 (2004), 1529--1550.


\end{thebibliography}
\end{document}